\documentclass[11pt]{amsart}%
\usepackage{amsfonts}
\usepackage{graphicx}
\usepackage{amscd}
\usepackage{amsmath}
\usepackage{amssymb}
\usepackage{amsfonts}
\newtheorem{theorem}{Theorem}[section]
\theoremstyle{plain}
\newtheorem{corollary}[theorem]{Corollary}
\newtheorem{definition}[theorem]{Definition}

\newtheorem{proposition}[theorem]{Proposition}

\theoremstyle{definition}
\newtheorem{example}[theorem]{Example}

\numberwithin{equation}{section}

\begin{document}

\title{Local regularity in non-linear generalized functions}

 \author[S. Pilipovi\'{c}]{Stevan Pilipovi\'c}
 \address{Department of Mathematics and Informatics, University of Novi Sad, Trg Dositeja Obradovi\'ca 4, Novi Sad, Serbia}
\email {stevan.pilipovic@dmi.uns.ac.rs}

  \author[D. Scarpalezos]{D. Scarpalezos}
\address{University Paris 7, Place Jussieu, Paris, France}
    \email{dim.scarpa@gmail.com}   

 \author[J. Vindas]{Jasson Vindas} 
\address{Department of Mathematics, Ghent University, Krijgslaan 281 Gebouw S22, B-9000 Gent, Belgium}
    \email{jvindas@cage.Ugent.be}       
       

 \begin{abstract}
In this review article we present regularity properties of generalized functions
which are useful in the analysis of non-linear problems. It is shown that  Schwartz distributions embedded
 into our new spaces of generalized functions, with additional properties described through the association,
 belong to various classical spaces with finite or infinite type of regularities.
\end{abstract}

\maketitle

%%%%%%%%%%%%%%%%%%%%%%%%%%%%%%%%%%%%%%%%%%%%%%%%%%%%%%%%%%%%%%%%%%
\section{Introduction}\label{sec0}
%%%%%%%%%%%%%%%%%%%%%%%%%%%%%%%%%%%%%%%%%%%%%%%%%%%%%%%%%%%%%%%%%%
Generalized function algebras of Colombeau type contain Schwartz's
distribution spaces and the embeddings preserve all the linear operations for
distributions. A great advantage of the generalized function algebra approach is
that various classes of nonlinear problems can be studied in these
frames as well as linear problems with different kinds of
singularities. We refer to \cite{bia},
\cite{col}, \cite{col1}, \cite{gkos} and \cite{ober001} for the theory of
generalized function algebras and their use in the study of
various classes of equations. For the purpose of local and
microlocal analysis, one is lead to study classical function spaces
within these algebras. 

In this survey paper we will present our investigations related to local analysis. We mention that in number of papers we have also studied 
wave front sets of  $\mathcal G^\infty$ and analytic types. The attention in article will be focused on regularity theory in  generalized function algebras. This regularity theory is parallel to the corresponding one within
distribution spaces related to analytic, real analytic, harmonic and  Besov type spaces, especially to  Zygmund type spaces.

Elements of  algebras of generalized functions
are represented by nets $(f_\varepsilon)_\varepsilon$ of smooth
functions, with appropriate growth as $\varepsilon\rightarrow 0,$ so that the spaces of Schwartz's distributions are embedded
into the corresponding algebras, and that for the space of smooth
functions the corresponding algebra of smooth generalized function
is $\mathcal G^\infty$ (see \cite{ober001}, \cite{Ver}). Elements of  these
algebras are obtained through the regularization of distributions
(convolving them with delta nets) and the construction of appropriate
algebras of moderate nets and null nets of smooth functions and
their quotients, as Colombeau did  (\cite{col1})  with his algebra $\mathcal G(\mathbb R^d)$, in such a way 
that distributions are included as well as their natural linear operations.

 The main goal of our investigations has been  to find out conditions with respect to the
growth order in $\varepsilon$ or integrability conditions with respect to $\varepsilon$
which characterize  generalized function spaces and algebras with finite type regularities.
 Our definitions for such generalized function spaces
 enable us to obtain information on the regularity properties of Schwartz 
 distributions embedded into the corresponding space of generalized functions. 
 
One can find many articles in the literature where local and
microlocal properties of generalized functions in generalized
function algebras have been considered. Besides the quoted monographs we
refer to the papers \cite{ar1}, \cite{ar2}, \cite{coga2}, \cite{coga}
\cite{gar1}--\cite{garhor}, \cite{SkalInvar}, \cite{horman}-\cite{OberGF00}, \cite{kuob}-\cite{psvv}, \cite{dim}, \cite{valm}. We also mention that some Tauberian theorems for regularizing transforms, \cite{pv}, \cite{p-r-v}, \cite{vipira}, \cite{DZ}, are also valuable tools for the study of regularity properties of generalized functions.

The paper is organized as follows. 
In Section \ref{gfc} we recall basic definitions of Colombeau type algebras, as well as basic notions related to Littlewood-Paley decompositions,
Besov and Zygmund type spaces. Holomorphic generalized functions are described in Section \ref{hgfa}, whereas the real analytic functions are presented in Section \ref{ragf}. We note that a fine contribution to the understanding of these classes has been given by Aragona and his collaborators, as well as Colombeau as a founder of the theory.  Differences between classical and generalized harmonic functions are presented as well as 
the fact that for the analytic and real analytic generalized functions standard points are sufficient for their description. This is not the case for the harmonic generalized functions, whose treatment needs additional preparation. 
Harmonic generalized functions with specific notions of $H-$ boundedness and generalized removable singularities are reviewed in Section
\ref{hgf}. Following \cite{bucu}, we discuss various properties of generalized functions with dependence on the  variable $\varepsilon$. Continuity, smoothness and the
  measurability condition with respect to $\varepsilon$ are discussed in Section \ref{classes}.  Actually in this section we introduce for the first time in the literature  the analysis of generalized functions 
through the integration with respect to $\varepsilon$. More precisely, we introduce in Section \ref{classes} the spaces $\mathcal G_q(\Omega) $, where the growth order is measure in an $L^{q}$ space with respect to the variable $\varepsilon$. In  Sections \ref{bes} and \ref{cass} we present Besov type spaces of generalized functions through the classes $G_q(\Omega). $ The regularity type results from Section \ref{cass} are rather recent. The last  section, Sections \ref{zcass}, is devoted to Zygmund type spaces. These spaces are essentially related to the process of regularizations  of Schwartz distributions which enable us to give precise characterizations of regularity properties of Schwartz distributions embedded into the corresponding space of Zygmund type generalized functions.  
  
   Let us note that many other  authors have made great contributions to developing the local and microlocal analysis of generalized functions of Colombeau type. Our collaborators Oberguggenberger, Vernaeve and Valmorin contributed very much to the results presented in this paper.  Among the main contributors to the local and micro local theory of Colombeau type generalized functions, we have to underline the role of
   Colombeau, Oberguggenberger, Aragona, H\" ormann, Kunzinger, Vernaeve, Garetto, Marti, Valmorin,  as well as their (and our ) coauthors.  In this sense, our list of references is rich enough but not complete.    
\section{Colombeau algebras and spaces}
\label{gfc}
Let $\Omega$ be an open subset of $\mathbb{R}^{d}$. We consider the families of local Sobolev seminorms
$||\phi||_{W^{m,p}(\omega)}= \sup
\{||\phi^{(\alpha)}||_{L^{p}(\omega)};\:|\alpha|\leq m \},$ where $m\in
{\mathbb N}_0$, $p\in[1,\infty]$, for $\omega\subset\subset\Omega$ (which means that $\overline{\omega}$ is compact in $\Omega$).

The spaces of moderate nets and negligible nets $\mathcal
E_{L_{loc}^p,M}(\Omega)$ and $\mathcal  N_{L_{loc}^p}(\Omega)$,
resp., $\mathcal E_{L^p,M}(\Omega)$ and $\mathcal  N_{L^p}(\Omega)$,
% and $\mathcalN_{L^p}(\Omega)$,
consist of nets $(f_{\varepsilon
})_{\varepsilon\in(0,1)}=(f_{\varepsilon
})_{\varepsilon}\in \mathcal{E}(\Omega)^{(0,1)}$ with the
properties (given by Landau's big $O$ and small $o$ )
\begin{equation}\label{2drs}
(\forall m\in\mathbb{N}_0)(\forall
\omega\subset\subset\Omega)(\exists a\in\mathbb{R})
(||f_{\varepsilon}||_{W^{m,p}(\omega)}=O(\varepsilon^{a}))
\end{equation}
\[
\mbox{ and } \;(\forall m\in\mathbb{N}_0)(\forall
\omega\subset\subset \Omega)(\forall b\in\mathbb{R})
(||f_{\varepsilon}||_{W^{m,p}(\omega)}=O(\varepsilon^{b})),
\]
The Sobolev lemma implies 
 $\mathcal
E_{L_{loc}^p}(\Omega)$ and $\mathcal  N_{L_{loc}^p}(\Omega)$ 
are algebras, but also that
$$\mathcal
\mathcal E_{M}(\Omega)=E_{L_{loc}^\infty}(\Omega)=\mathcal
E_{L_{loc}^p}(\Omega),\;\; \mathcal
N(\Omega)=\mathcal
N_{L_{loc}^\infty}(\Omega)=\mathcal N_{L_{loc}^p}(\Omega),\ p\geq 1.$$ Thus
the Colombeau algebra of generalized functions can be defined as
$$
\mathcal
G(\Omega)= \mathcal E_{L_{loc}^p}(\Omega)/\mathcal N_{L_{loc}^p}(\Omega), \  \ p\geq1.
$$
Recall \cite{ober001}, that the algebra of regular generalized functions as defined as
$\mathcal{G}^\infty(\Omega)=
\mathcal{E}^\infty_{M}(\Omega)/\mathcal{N}(\Omega),$ where
$\mathcal{E}^\infty_{M}(\Omega)$ consists of nets $(f_{\varepsilon
})_{\varepsilon\in(0,1)}\in \mathcal{E}(\Omega)^{(0,1)}$ with the
property
$$
(\forall K\subset\subset \Omega) (\exists a\in\mathbb{R}) (\forall
n\in\mathbb{N}) (|\sup_{x\in
K}f^{(n)}_{\varepsilon}(x)|=O(\varepsilon^{a})).
$$

If the elements of the nets $(f_\varepsilon)_\varepsilon\in
\mathcal E_{M}(\Omega)$  are constant functions in $\Omega$ (i.e.,
seminorms  reduce to the absolute value), then one obtains the
corresponding algebras $\mathcal{E}_0$ and $\mathcal{N}_0;$
$\mathcal{N}_0$ is an ideal in $\mathcal{E}_0$ and, as their
quotient, one obtains the Colombeau algebra of generalized complex
numbers: $\bar{\mathbb{C}}=\mathcal{E}_0/\mathcal{N}_0$ (or
$\bar{\mathbb{R}}$ if nets are real). It is a ring, not a field.
For the analysis of harmonic generalized functions we will recall the definition of compactly supported generalized points, due to
Oberguggenberger and  Kunzinger in \cite{kuob}.
A net $(x_\varepsilon)_\varepsilon$ in a general metric space $(A,d)$ is called moderate, if
\[
(\exists N\in\mathbb N)(\exists x\in
A)(d(x,x_\varepsilon)=O(\varepsilon^{-N})).
\]
and an equivalence relation in $A^{(0,1]}$ is introduced by
\[
(x_\varepsilon)_\varepsilon\sim
(y_\varepsilon)_\varepsilon\Leftrightarrow (\forall p\geq 0)
(d(x_\varepsilon,y_\varepsilon)=O(\varepsilon^p)).
\]
$\widetilde A=A/\sim$ is called the set of generalized points in
$A$.
 If $A=\Omega$ is an opens subset of
$\mathbb R^d$, then $\widetilde\Omega=\Omega/\sim$ is the set of
generalized points.  Note $\tilde{\mathbb{C}}=\bar{\mathbb{C}}$
($\tilde{\mathbb{R}}=\bar{\mathbb{R}}$).
An element $\widetilde x\in\widetilde{\Omega}$
is called compactly supported if $x_\varepsilon$ lies in a compact
set for
 $\varepsilon<\varepsilon_0$ for some $\varepsilon_0\in(0,1)$. The set of compactly supported
points $\widetilde x $ ($\in\widetilde{\Omega}$) is denoted by
$\widetilde{\Omega}_c$.  Recall, nearly
standard points are elements $\widetilde x\in\widetilde{\Omega}_c$
with limit in $\Omega$ that is, there exists $x\in\Omega$ such
that for a representative  $(x_\varepsilon)_\varepsilon$,
$x_\varepsilon\rightarrow x$, as $\varepsilon\rightarrow 0$,
holds. We denote by $\widetilde{\Omega}_{ns}$ the set of nearly
standard points of $\Omega.$

The embedding of the Schwartz distribution space
$\mathcal{E}^{\prime}(\Omega)$ into $\mathcal{G}(\Omega)$ is realized through the sheaf
homomorphism $ \mathcal{E}^{\prime}(\Omega)\ni T\mapsto\iota(T)=
[(T\ast\phi_{\varepsilon}|{\Omega})_\varepsilon]\in\mathcal{G}(\Omega)$,
where the fixed net of mollifiers
$(\phi_{\varepsilon})_{\varepsilon}$ is defined by
$\phi_{\varepsilon}=\varepsilon ^{-d}\phi(\cdot/\varepsilon),\;
\varepsilon<1,$ and $\phi\in\mathcal{S}(\mathbb{R}^{d})$ satisfies
$$\int_{\mathbb{R}^{d}}
\phi(t)\mathrm{d}t=1,\;\int_{\mathbb{R}^{d}}
t^{m}\phi(t)\mathrm{d}t=0,\: |m|>0.$$
($t^{m}=t_{1}^{m_{1}}...t_{d}^{m_{d}}$ and $|m|=m_{1}+...+m_{d}$).

This sheaf homomorphism \cite{gkos}, extended over $\mathcal{D}^{\prime}$,
gives the embedding of $\mathcal{D}^{\prime}(\Omega)$ into
$\mathcal{G}(\Omega)$. 
We also use the notation $\iota$ for the mapping from $\mathcal{E}'(\Omega)$
into $ \mathcal{E}_M(\Omega)$, $\iota(T)=
(T\ast\phi_{\varepsilon}|{\Omega})_\varepsilon$. Throughout this article, 
$\phi$ is fixed and satisfies the above condition over its moments.

We will use a continuous Littlewood-Paley decomposition of the unity
(see \cite{triebel2006} or \cite[Sect. 8.4]{hore}, for instance).
Let $\varphi\in\mathcal S(\mathbb R^d)$ such that its Fourier transformation
$\hat{\varphi}\in
{\mathcal D}({\mathbb R^{d}})$ is a real valued radial
(independent of rotations) function with support contained in
the unit ball such that $\hat{\varphi}(y)=1$ if $|y|\leq 1/2$. Set
$\displaystyle \hat{\psi}(y)=-\frac{d}{d\varepsilon}\hat{\varphi}(\varepsilon
y)|_{\varepsilon=1} =-y\cdot \nabla\hat{\varphi}(y).$ The support of
$\hat{\psi}$ is contained in the set $1/2\leq |y|\leq1$.
Then,
$\varphi$ has the same moment properties as $\phi$ given above but in the sequel we will use term mollifier only for $\phi$
since it will be used for the embedding of distributions into various spaces and algebras of generalized functions. The function
  $\psi$ is a wavelet (all the moments of $\psi$ are equal to zero). 
  
One has \cite[Sect. 8.6]{hore}: for any $u\in {\mathcal S}'(\mathbb R^d)$,
\begin{equation}\label{decompose11}
u=u*\varphi+\int_{0}^{1}u*\psi_\eta \frac{d\eta}{\eta}=u*\varphi_{\varepsilon}+\int_0^{\varepsilon}u*\psi_\eta \frac{d\eta}{\eta},\
0<\varepsilon\leq 1,
\end{equation}
and hence,
\begin{equation}\label{decompose}
u\ast \varphi_\varepsilon=u*\varphi+\int_{\varepsilon}^{1}u*\psi_\eta \frac{d\eta}{\eta},\
0<\varepsilon\leq 1.
\end{equation}
We recall that the  Besov spaces $
B^{s}_{q,p}(\mathbb R^d), p,q\in(0,\infty], s\in \mathbb R,$
are defined as
$$B^{s}_{q,p}(\mathbb R^d)=\{f\in\mathcal S'(\mathbb R^d):
||u||^s_{q,p}
<\infty\}, \mbox{ where }$$
\begin{equation}\label{besov} 
||u||^s_{q,p}: = ||f*\phi||_{L^p(\mathbb R^d)}+ \left(\int_0^1y^{-sq}||f*\psi_y||^q_{^{L^p(\mathbb R^d)}}
dy/y\right)^{1/q}.
\end{equation}
The definition is independent of the choice of the pair $\varphi$. If we discretize $y=2^{-j}, j\in \mathbb N,$ then one can replace \cite{triebel2006} the second term in (\ref{besov})
by
 $$\left(\sum_{j=1}^\infty2^{-jqs}||f\ast\psi_{2^{-j}}||_{L^p(\mathbb R^d)}^q\right)^{1/q}.
 $$

They are quasi-Banach spaces if $\min\{p,q\}\leq 1,$ and Banach spaces if $\min\{p,q\}\geq 1.$ In the sequel we consider the second case i.e.   $\min\{p,q\}\geq 1.$

In particular, the Zygmund spaces are defined by
$$%\begin{equation}\label{hz}
C_{\ast}^{s}(\mathbb R^d):=\{u\in {\mathcal S}':\;
|u|_{\ast}^{s}:=||\phi*u||_{L^{\infty}}
+\sup_{y\in(0,1)}(y^{-s}||\psi_{y}*u||_{L^{\infty}})<\infty\},
$$%\end{equation}
so that $C_{\ast}^{s}(\mathbb R^d)=B^{s}_{\infty,\infty}(\mathbb{R}^{d})$.
We will also  consider 
spaces $B^{s}_{q,L^p_{loc}}(\Omega)$ consisting of tempered distributions with the property:
$f\in B^{s}_{q,L^p_{loc}}(\Omega)$ if $\theta f\in B^{s}_{q,p}(\mathbb R^d)$ for every $\theta\in \mathcal{D}(\Omega).$

\section{Holomorphic generalized functions \cite{pili}}\label{hgfa}

We denote by ${\mathcal O}(\Omega)$ the space of holomorphic functions
on $\Omega$, where $\Omega$ is an open subset of
$\mathbb{R}^2=\mathbb{C}$; $D(z_0,r)$ denotes a disc with the
center $z_0$ and radius $r>0.$
\begin{definition}
A generalized function
$f=[(f_\varepsilon)_\varepsilon]\in{\mathcal G}(\Omega)$ is said to be holomorphic if $\frac{\partial f}{\partial %%@
\bar{z}}=0$ in ${\mathcal G}(\Omega)$.

The set of holomorphic generalized functions is denoted by ${\mathcal
G}_H(\Omega)$. 
\end{definition}
Due to
Colombeau-Gal\'{e} \cite{coga2}
we know that
$f\in{\mathcal G}_H(\Omega)$ if and only if for every relatively
compact open set $\Omega'$ in $\Omega$, $f$ admits a
representative $(f_{\varepsilon})_{\varepsilon} \in{\mathcal
E}_M(\Omega')$ with $f_{\varepsilon}\in{\mathcal O}(\Omega')$,
$\varepsilon\in(0,1]$.
In the same paper it is shown that ${\mathcal{G}}_{H}(\Omega) \cap
\mathcal{D}'(\Omega)={\mathcal{O}}(\Omega)$.
Holomorphic generalized functions can be well understood  through the next theorem.
\begin{theorem} 
Let $(f_{\varepsilon})_{\varepsilon}\in{\mathcal{E}}(\Omega)^{(0,1]}$ and suppose
that for every point $z_0\in\Omega$ there exist $r_\varepsilon>0,
0<\varepsilon\le 1$ such that
$$
\begin{array}{rl}
(i) & \bar{\partial}{f_\varepsilon}|_{D(z_0,r_\varepsilon)}=0, \;
0<\varepsilon\le 1, \mbox{i.e. the restrictions to the discs vanish}.
\\
(ii) & \exists \eta>0, \exists a>0,  \exists \varepsilon_0\in(0,1],
\,\,|f_{\varepsilon}^{(n)}(z_0)|\le \eta^{n+1}n!\varepsilon^{-a},\,
n\in\mathbb{N}, \varepsilon\in(0,\varepsilon_0).
\end{array}
$$
Then $(f_\varepsilon)_\varepsilon\in{\mathcal{E}_M}(\Omega)$ and
$[(f_\varepsilon)_\varepsilon] \in{\mathcal{G}_H}(\Omega)$.
\end{theorem}

This leads to the  result that  ${\mathcal
G}_H(\Omega)\subset{\mathcal G}^\infty(\Omega),$ where the algebra of regular generalized functions ${\mathcal G}^\infty(\Omega)$ was defined in Section \ref{gfc}

%
%this mapping is extended on ${\cal D}'(\Omega).$
%Note that $\mathcal{G}^\infty$ is a subsheaf of $\mathcal{G}.$

%(see \cite{ober001} for the definition of ${\mathcal
%G}^\infty(\Omega).$ 
%

We collect properties of holomorphic generalized functions in the
next theorem.

\par 

\begin{theorem} 
Let
%$\Omega$ denote an open set in $\mathbb{R}^2$ and
$g\in{\mathcal G}_H(\Omega)$.

(i) $g$ admits a representative $(g_{\varepsilon})_{\varepsilon}$ such that
$g_{\varepsilon}\in{\mathcal O}(\Omega),$ $\varepsilon\in(0,1]$.

(ii) $g=0$ if and only if for any open set $\Omega'\subset\subset
\Omega$ there exists a representative $(g_\varepsilon)_\varepsilon\in {\mathcal
O}(\Omega')^{(0,1)}$ such that
\begin{equation} \label{vt}
\begin{array}{l}
\forall z_0\in\Omega',\forall a>0,\exists \eta>0,   \exists
\varepsilon_0\in(0,1],
\exists C>0, \vspace{4pt}\\
|g_{\varepsilon}^{(n)}(z_0)|\le \eta^ {n+1}n!\varepsilon^{a},\,\,
n\in\mathbb{N}, \varepsilon\in(0,\varepsilon_0).
\end{array}
\end{equation}

(iii) Let $z_0\in\Omega$, $\eta > 0$. Then $f$ vanishes in the
disc $V = D(z_0,\frac{1}{\eta})$ if and only if there exists a
representative $(g_\varepsilon)_\varepsilon$ of $f|V$ in ${\mathcal
O}(V)^{(0,1]}$ such that
$$
\forall a>0,  \exists \varepsilon_0\in(0,1] :
|g_{\varepsilon}^{(n)}(z_0)|\le \eta^ {n+1}n!\varepsilon^{a};\,\,
n\in\mathbb{N}, \varepsilon\in(0,\varepsilon_0).
$$
\end{theorem}

%\begin{corollary}
%Let $f\in{\cal G}_H(\Omega)$. Then $f$ vanishes in a non-void open subset of $\Omega$ if and only if there exists  %%@
%$z_0\in\Omega$, an open neighborhood $V_0$ of $z_0$ and a representative $(g_\varep)_\varep$ of $f|V$ in ${\cal %%@
%O}(V)^{(0,1]}$ such that
%$$
%(**)\,\,\,\exists \eta>0, \forall a>0,  \exists \varep_0\in(0,1),
%|g_{\varep}^{(n)}(z_0)|\le \eta^ {n+1}n!\varep^{a};\,\,
%n\in\mathbb{N}, \varep\in(0,\varep_0).
%$$
%\end{corollary}

As a consequence, we obtain  a simple proof  that on a connected open set  
$\Omega$ an $f\in{\mathcal G}_H(\Omega)$ vanishes on $\Omega$ if and only if it
vanishes on a non-void open subset of $\Omega$. This is an important result of Khelif and Scarpalezos \cite{kelsca}. It shows a main property of holomorphic generalized functions.
For their analysis, standard points are enough while for
Colombeau type generalized functions, generalized points are essential \cite{kuob} (see also \cite{KKo1} for more general aspects of this fact.)

 The existence of a global holomorphic representative of $f\in\mathcal{G}_H(\Omega)$
depends on $\Omega\subset \mathbb{C}^d$ is still an open problem although
such representation for appropriate domains
$\Omega\subset \mathbb{C}^n$ can be constructed. In the one dimensional case we have
\begin{theorem}\label{ghr}
If $f\in \mathcal G_H(\mathbb C),$ then there exists one of its representative $(f_\varepsilon)_\varepsilon$
consisting of entire functions on $\mathbb C.$
\end{theorem}

%\section{Global representatives for solutions of elliptic equations}

\section{Real analytic generalized functions \cite{psv}} \label{ragf}

Now we are considering  $\omega$, an open set in ${\mathbb R^d}$.
\begin{definition}\label{defi, aaa}
Let $x_0\in \omega$. A  generalized function $f\in{\mathcal G}(\omega)$ is said to be real
analytic at $x_0$ if there exist an open ball $B=B(x_0,r)$ in $\omega$
containing $x_0$ and $(g_\varepsilon)_\varepsilon\in {\mathcal E}_M(B)$ such that
$$
(i)\,\,\,\,\,\,\, f|B=[(g_\varepsilon)_\varepsilon]\,\,{\rm in}\,\,{\mathcal G}(B);
$$
$$
(ii)\,\,\,\,\,\,(\exists \eta>0)(\exists a>0)(\exists \varepsilon_0\in
(0,1))
$$
$$
\sup_{x\in {B}}|\partial^\alpha g_\varepsilon(x)| \le
\eta^{|\alpha|+1}\alpha!\varepsilon^{-a},\; 0<\varepsilon<\varepsilon_0,\; \alpha\in\mathbb N^d.
$$
It is said that $f$ is real analytic in $\omega$ if $f$ is real
analytic at each point of $\omega$.
The space of all generalized functions which are real analytic in $\omega$ is denoted by %%@
$\mathcal{G}_A(\omega)$.

The analytic singular support, $\operatorname*{singsupp }_{ga} f,$ is the
complement of the set of points $x\in \omega$ where $f$ is real
analytic.
\end{definition}

\noindent It follows  from the definition that ${\mathcal G_A}$ is a subsheaf of ${\mathcal %%@
G}$.

Using Stirling's formula it is seen that condition (ii) in Definition \ref{defi, aaa} is %%@
equivalent to
$$
(iii)\,\,\,\,\,\,
(\exists \eta>0)(\exists a>0)(\exists \varepsilon_0\in (0,1))
$$
$$
\sup_{x\in B}|\partial^\alpha g_\varepsilon(x_0)| \le
\eta(\eta|\alpha|)^{|\alpha|}\varepsilon^{-a},\; 0<\varepsilon<\varepsilon_0,\;
\alpha\in\mathbb N_0^d.
$$
The use of Taylor expansion and condition $(ii)$ of Definition
\ref{defi, aaa} imply that $(g_{\varepsilon})_\varepsilon$ admits a
holomorphic extension in a complex ball $B=\{z\in
\mathbb{C}^d; |z-x_0|<r \}$ which is independent of $\varepsilon$.
Consequently we get a holomorphic extension $G$ of $[(g_\varepsilon)_\varepsilon]$
and then $f|B=G|B$. (It is clear from the context whether $B$ is a complex or real ball.)

 The existence of a global real analytic representative of $f\in\mathcal{G}_A(\mathbb R^d)$, as in the case of analytic generalized functions
 is also an open problem. In the case $d=1$ we have a positive answer, there exist a global representative, i.e. a representative a
 real analytic generalized function on $\mathbb R$ consisting of real analytic functions defined on $\mathbb R$.

Moreover, we have a similar situation for real analytic generalized functions as for holomorphic, they are determined by their values at standard points.
We have
\begin{theorem}

a) Let $\Omega$ be an open set of $\mathbb{C}^p, p>1$ and $f=[(f_\varepsilon)_\varepsilon]\in \mathcal G_H(\Omega)$ such that
$f(x)=0$ for every  $x\in \Omega$ ($(f_\varepsilon(x))_\varepsilon\in\mathcal{N}(\Omega)$). Then $f\equiv 0.$

b) Let $\omega$ be an open set of $\mathbb{R}^d, d\in\mathbb{N}$
and $f=[(f_\varepsilon)_\varepsilon]\in \mathcal G_A(\omega)$ such that
$f(x)=0$ for every  $x\in \omega$. Then $f\equiv 0.$
\end{theorem}

The singular support of an $f\in \mathcal G(\omega)$ is defined as the complement of the union of open sets in $\omega$ where $f\in\mathcal G^\infty(\omega).$ In a similar way one defines the notion
$\mbox{ singsupp }_{ga}f.$ 
\begin{theorem} 
\cite{psv}
Let $f\in\mathcal{E}'(\omega)$ and $f_\varepsilon=f*\phi_{\varepsilon}, \varepsilon\in (0,1)$
be its regularization by a net $\phi_\varepsilon=\varepsilon^{-1}\phi(\cdot/\varepsilon), \varepsilon<1,$
where $(\phi_\varepsilon)_\varepsilon$ is a net of mollifiers. Then
$$\operatorname*{singsupp }_a f = \operatorname*{singsupp }_{ga} [(f_\varepsilon)_\varepsilon].
$$
%where on the left hand side is the analytic singular support of the distribution $f$.
\end{theorem}

\section{Harmonic generalized function \cite{piscahr}}\label{hgf}

%Our results of this section  are proved in  \cite{piscahr}.

We denote by $Har(\Omega)$ the space of harmonic functions in
$\Omega$.
\begin{definition} \label{dh}
We call a generalized function $G\in\mathcal G(\Omega)$ harmonic
generalized function, if $\Delta G=0$ holds in $\;\mathcal
G(\Omega)$. The linear space of harmonic generalized functions in
$\Omega$ is denoted by $\mathcal{G}_{Har}(\Omega)$.
\end{definition}
We have shown that  $\mathcal{G}_{Har}$
is a closed subsheaf of the sheaf of
$\widetilde{\mathbb{C}}$-modules $\mathcal{G}$.

Moreover we have the next important result
\begin{theorem}
Every harmonic generalized function $G\in \mathcal{G}(\Omega)$
admits a global harmonic representative $(G_\varepsilon)_\varepsilon$,
that is, for each $\varepsilon\leq 1$, $G_\varepsilon$ is
harmonic.
\end{theorem}
We call the $(G_\varepsilon)_\varepsilon$ {a global harmonic
representative} of $G$. For the main results of our analysis we use global harmonic representatives.

Obviously,
 $\Omega\rightarrow \mathcal{G}_A(\Omega)$ is a  subsheaf of the sheaf of algebras
$\mathcal{G}^\infty.$
Moreover, every harmonic generalized function is a real analytic generalized
function.
Since $\mathcal G_{Har}(\Omega)$ is a submodule of
$\mathcal G_A(\Omega)$ the following consequence  is immediate
(see the previous section).
\begin{theorem}
\label{prv} Let $\Omega $ be a connected open subset of ${\mathbb
R}^d$ and $f \in \mathcal G_{Har}(\Omega)$.  If there exists
$A\subset \Omega$ of positive Lebesgue measure ($\mu(A)>0$) such
that $f(x)=0$ for every $x\in A$, then $f\equiv 0$.
\end{theorem}

For the generalized maximum principle, we need additional notation. 
Let $K$ be a compact set of $\Omega$, $\widetilde{x}_0\in
\widetilde{\Omega}_c$ be supported by $K$  and  let $r>0$ such
that $K+B(0,r)\subset \subset \Omega. $ With such an
$\widetilde{x}_0$, we denote by $B(\widetilde{x}_0,r)$ the
following subset of $\widetilde{\Omega}_c$
\begin{equation} \label{kolo}
B(\widetilde{x}_0,r)=\{\widetilde{t}=[(t_\varepsilon)_\varepsilon]\in\widetilde{\Omega}_c;
|x_{0,\varepsilon}-t_\varepsilon|\leq r,\; \varepsilon\leq 1\}.
\end{equation}
We call the set
$B(\widetilde{x}_0,r)\subseteq\widetilde{\Omega}_c$ a semi-ball in
$\Omega$.  Note that by now we distinguish between balls
$B(x_0,r)\subset \mathbb K$, balls $\widetilde B(\widetilde
x_0,r)$ in $\widetilde{\mathbb K}$ and semi-balls $B(\widetilde
x_0,r)$.
%Note that $S_{-\rho}$ is an open connected set.
\begin{theorem} [Maximum principle]\label{maxp} 
\label{max1}    Let $G$ be  a real-valued harmonic generalized
function in an open set $\Omega$. Then the following holds:
\begin{enumerate}
\item \label{constant1} Let $r>0$ and $\widetilde{x}_0\in\widetilde{\Omega}_c$
with representative $(x_{0,\varepsilon})_\varepsilon$ be given
such that  $B({x}_{0,\varepsilon},2r) \subset\subset \Omega$,
$\varepsilon\leq 1$ (that is the semi-ball
$B(\widetilde{x}_0,2r)\subset \widetilde{\Omega}_c$). Suppose
$G(\widetilde{x}_0)\geq G(\widetilde{t})$ for each
$\widetilde{t}\in B(\widetilde{x}_0,r)$. Then $G$ is a constant
generalized function in $B(\widetilde{x}_0,r)$.
 \item \label{constant2} Let
$\Omega$ be connected. If there exists a compactly supported point
$\widetilde{x}_0\in\widetilde{\Omega}_c$ such that
 $G(\widetilde{x}_0)\geq
G(\widetilde{t})$, $\widetilde{t}\in \widetilde{\Omega}_c$, then
$G$ is a constant generalized function in $\Omega$.
\end{enumerate}
\end{theorem}

\begin{proposition}\label{invm}
Let $G=[(G_\varepsilon)_\varepsilon]\in \mathcal G(\Omega)$ such that for every
$\widetilde{x}=[(x_\varepsilon)_\varepsilon] \in \widetilde{\Omega}_c$ and every
$R>0$ such that the semi-ball $\tilde{B}(\widetilde{x},R)\subset
\widetilde{\Omega}_c,$
\begin{equation} \label{nrep}
[(G_\varepsilon(\widetilde{x}))_\varepsilon]=[(\frac{1}{V_R}\int_{B(x_\varepsilon,R)}
G_\varepsilon(t)dt)_\varepsilon].%\; \eps\leq 1,
\end{equation}
%where $(x_\eps)_\eps$ is a representative of $\widetilde{x}$. ((\ref{nrep}) does not depend on a representative.)
Then $G\in \mathcal G_{Har}(\Omega). $
\end{proposition}

\begin{theorem}\label{um} Let $\Omega$ be  connected.

(i) With the assumptions of Theorem \ref{maxp} (i),
 $G$ is a constant generalized function in
$\Omega.$

(ii) If there exists a nearly standard point
$\widetilde{x}_0\in\widetilde{\Omega}_{ns}$ such that
$G(\widetilde{x}_0)\geq G(\widetilde{t})$ for each nearly standard
point $\widetilde{t}\in \widetilde{\Omega}_{ns}$, then $G$ is a
constant generalized function in $\Omega$.
\end{theorem}

\begin{corollary}
(i) Let $u$ be a complex harmonic generalized function in  a connected
open set $\Omega$. If $|u|$ has a maximum
$\widetilde{M}\in\widetilde{\mathbb R}$  at  $\widetilde{x}_0\in
\widetilde{\Omega}_c$, then $u\equiv
\widetilde{A}=u(\widetilde{x}_0)\in \widetilde{\mathbb C}$,
$|\widetilde{A}|=\widetilde{M}$.
%The following assertion is a consequence of the previous one.
%\begin{corollary} \label{sub1}

(ii) Let $G\in{\mathcal G}_{Har}(\Omega)$ be a non-constant real valued
generalized function. Then
\begin{enumerate}
\item \label{maxim1} $G$ does not attain its maximum inside $\Omega$, that is at
a generalized point $\widetilde t_0\in\widetilde{\Omega}_c$. \item
\label{maxim2} Let $\Omega'\subset\subset\Omega$ be open. Then the
maximum of $G$ in $\overline{\Omega'}$ is attained at a
generalized point supported by the boundary of $\Omega'$.
\end{enumerate}
\end{corollary}
For the generalizations of Liouville's theorem for harmonic generalized 
functions we need to repeat several notions for  $u\in\mathcal G(\mathbb R^d)$. 

We call $u$:

(i) non-negative, if for each compact set $K$ there exists a
representative $(u_\varepsilon)_\varepsilon$ of $u$ such that for each $\varepsilon>0$,
$\inf_{x\in K}u_\varepsilon(x)\geq 0$.

(ii) strictly positive, if for
each representative $(u_\varepsilon)_\varepsilon$ of $u$ and for each compact
set $K$ there exists constants $m$ and  $\varepsilon_0$ such that for
each $\varepsilon<\varepsilon_0$, $\inf_{x\in K}u_\varepsilon(x)\geq \varepsilon^m$.

(iii) A harmonic generalized function $u$ is
said to be globally non-negative, if it admits a global harmonic
representative $(G_\varepsilon)_\varepsilon$ so that $G_\varepsilon$ is non-negative
for each $\varepsilon\leq 1$.

(iv) We
call a harmonic generalized function $u$ H-non-negative  (and
write $u\geq_H 0$), if it admits a global harmonic representative
$(G_\varepsilon)_\varepsilon$ with the following property:
\begin{equation}\label{Ppositive}
(\forall m>0)(\forall a>0)(\exists \varepsilon_{a,m}\in(0,1])(\forall
\varepsilon<\varepsilon_{a,m})( \forall\, t: \vert t \vert < \frac
{1}{\varepsilon^m})(u_\varepsilon(t)+\varepsilon^a\geq 0)
\end{equation}
Furthermore, a harmonic generalized function $u$ is said to be
H-bounded from above (resp.\ below)  by $\widetilde{c}\in
\widetilde{\mathbb R}$, if for a representative $(c_\varepsilon)_\varepsilon $
of $\widetilde{c}$, the global harmonic representative
$(G_\varepsilon)_\varepsilon-(c_\varepsilon)_\varepsilon$ satisfies condition
(\ref{Ppositive}).
%such that for each compact set $K$ there exists an $\eps_0$ such
%that for each $\eps<\eps_0$ we have $\inf_{x\in K}u_\eps(x)\geq 0$
%(resp.\ $\sup_{x\in K}u_\eps(x)\leq 0$).
%\item globally strictly positive, if there exists a representative $(u_\eps)_\eps$ of $u$ such that for each $\eps>0$
%and for all $x\in \mathbb R^n$ we have $u_\eps(x)\geq 0$.
A harmonic generalized function $u$ is said to be H-bounded if it
is H-bounded from above and from below.

\begin{theorem}\label{thpilliouville}
A harmonic generalized  function $u$ in $\mathbb R ^d$ which is
H-bounded from below is a constant.
\end{theorem}

A direct consequence  is that
every H-bounded harmonic generalized function $u\in {\mathcal
G}({\bf R}^d)$ is a constant.
%\end{corollary}

Now we give the definition of
isolated singularity of harmonic generalized function.
%\begin{definition} 

Let $\Omega$ be an open set of $\mathbb{R}^d$ and
$x_0\in\Omega.$
 A generalized function  $G\in\mathcal{G}(\Omega\setminus\{x_0)\}$
 (resp.\ $G\in\mathcal{G}_{Har}(\Omega\setminus\{x_0)\}$)
is said to have an isolated (resp.\ isolated harmonic) singularity
at $x_0$. Moreover, if there exists $F\in\mathcal{G}(\Omega)$
(resp.\ $F\in\mathcal{G}_{Har}(\Omega)$) such that
 $F|_{\Omega\setminus\{x_0\}}=G$, then it is said that $G$ has a removable (resp.\ harmonic removable) singularity.
%\end{definition}

Theorem \ref{hi} below states assertions on harmonic generalized
functions in pierced domains. First we need a definition of
$H$-boundedness in a neighborhood of $x_0$ which corresponds to a
$H$-boundedness at infinity.
%\begin{definition}\label{hzero}

Let $G\in \mathcal G_{Har}(\Omega\setminus\{x_0\})$ and let
$B(x_0,R)\subset \Omega.$  It is said that it is $H-$bounded in a
neighborhood of $x_0$ if there exists
$M=[(M_\varepsilon)_\varepsilon]>0$
 and a global harmonic
representative $(G_\varepsilon)_\varepsilon$ in
$\Omega\setminus\{x_0\}$ such that for every $m\in \mathbb{N}$
there exists $\varepsilon_m\in (0,1]$ such that
$$|G_\varepsilon(x)|<M_\varepsilon, x\in \{\varepsilon^m<|x-x_0|<R,
\varepsilon<\varepsilon_m\}.
$$

\begin{theorem}\label{hi}
Let $G\in\mathcal{G}_{Har}(\Omega\setminus\{x_0\})$. The following
holds:
\begin{enumerate}
\item \label{singer1} Assume additionally that $G\in\mathcal G(\Omega)$, and
that for every sharp neighborhood $V$ of $x_0$ $G$ has a
representative $(G_\varepsilon)_\varepsilon$ so that for every
$\varepsilon\leq 1,$ $G_\varepsilon $  is harmonic outside
$V_\varepsilon,$ where $V=[(V_\varepsilon)_\varepsilon]$. Then
$G\in\mathcal{G}_{Har}(\Omega)$.
\item \label{singer2} If $G$ is $H-$bounded at $x_0$,
then $G$ extends uniquely to an element of
$\mathcal{G}_{Har}(\Omega)$.
\end{enumerate}
\end{theorem}

\section{New spaces defined by integration in $\varepsilon$ \cite{psvv} }
\label{classes}

%\subsection{Measurable dependence on $\varepsilon$}
Following \cite{bucu}, we will consider representatives $(f_\varepsilon)_\varepsilon$, $(\varepsilon,x)\mapsto f_\varepsilon(x)$ 
which continuously depend on $\varepsilon$ or (moreover) smoothly depend on $\varepsilon\in (0,1]$
(always smooth in $x$). The notation $_{co}$ stands for the continuous parametrization, while 
$_{sm}$ stands for the smooth parametrization. It is obvious that
$$\mathcal E_{M,sm}(\Omega)\subset\mathcal E_{M,co}(\Omega)\subset\mathcal E_{M}(\Omega),
\mathcal N_{sm}(\Omega)\subset\mathcal N_{co}(\Omega)\subset\mathcal N(\Omega),
$$
Furthermore, it is shown in \cite{bucu} that
$$\mathcal G_{co}(\Omega)=\mathcal G_{sm}(\Omega)\subset\mathcal G(\Omega),
$$
where the last  inclusion is strict. The same relations hold for generalized complex (and real)
numbers.

Moreover, we will consider representatives $(f_\varepsilon)_\varepsilon$ which are measurable functions with respect to
 $\varepsilon:$  
$$
\mbox{for every fixed } x\in\Omega, (0,1)\ni \varepsilon\mapsto f_\varepsilon(x)\in\mathbb C \mbox{ is measurable}.
$$
%We use subscript $_{me}$ for the corresponding algebras and  have the next assertion.

Let $p\in[1,\infty].$ The definitions of algebras $\mathcal E_{L^p_{loc},M}(\Omega)=\mathcal E_M(\Omega)$ and $\mathcal N_{L^p_{loc}}(\Omega)
=\mathcal N(\Omega)$ from Section \ref{gfc} can be formulated by measurable representatives  with respect to $\varepsilon$. 
We will denote them by  the symbols
 $\mathcal E_{M,me}(\Omega), \mathcal N_{me}(\Omega)$ and their quotient by $\mathcal G_{me}(\Omega)$, where we assume measurability dependence. The next proposition is also from \cite{bucu}.

%\marginpar{(i) was modified}
\begin{proposition}\label{spec}
The following strict embeddings hold
$$\mathcal E_{M,co}(\Omega)\subset\mathcal E_{M,me}(\Omega);\; 
\mathcal N_{co}(\Omega)\subset\mathcal N_{me}(\Omega);\;
\mathcal G_{co}(\Omega)\subset\mathcal G_{me}(\Omega)\subset\mathcal{G}(\Omega).$$
\end{proposition}

 Let $p\in[1,\infty],$ $(f_\varepsilon)_\varepsilon\in\mathcal E_{M,me}(\Omega)$ and $\omega\subset\subset\Omega$. Then it is clear that
$(0,1)\ni\varepsilon\mapsto ||f_\varepsilon(\cdot)||_{L^p(\omega)}\in\mathbb R
$
is a measurable function.

%\marginpar{part (ii) can simply put as a comment after the proposition}
\begin{example}
Let $A\subset(0,1]$ be the well known non-measurable Vitaly set. Let $f_\varepsilon=1,\varepsilon\in A$, 
$f_\varepsilon=0, \varepsilon\in(0,1]\setminus A$. Then $[f_\varepsilon]$
 shows that  $\mathcal G_{me}(\Omega)$ is strictly contained in $\mathcal G(\Omega).$
\end{example}

Let $q\in[1,\infty)$ and  $p\in[1,\infty]$. We  say that a net of 
 $(f_\varepsilon)\in{\mathcal E}(\Omega)^{(0,1)}$ belongs to $\mathcal E_{q,L^p_{loc}}(\Omega)$, 
 respectively to $\mathcal N_{q,L^p_{loc}}(\Omega)$ if it is measurable, locally bounded on $(0,1],$ with respect to $\varepsilon$,
 for every fixed $x\in\Omega$ and it satisfies the growth estimates 
%\marginpar{there was a mistake, the moderateness subindex was not needed ``$M$''}
\begin{equation}
\label{eql1}
(\forall k\in\mathbb N_0)(\forall \omega\subset \subset \Omega)(\exists s\in \mathbb R)
(\int_0^1\varepsilon^{qs}||f_\varepsilon||_{W^{k,p}(\omega)}^qd\varepsilon/\varepsilon<\infty),
\end{equation}
%\marginpar{modified}
respectively,
\begin{equation}
\label{eql2}
(\forall k\in\mathbb N_0)(\forall \omega\subset \subset \Omega)(\forall s\in \mathbb R)
(\int_0^1\varepsilon^{qs}||f_\varepsilon||_{W^{k,p}(\omega)}^qd\varepsilon/\varepsilon<\infty).
\end{equation}
By the Sobolev lemma, it follows that these spaces of nets are independent of the value of $p$. We therefore set
$$
\mathcal{E}_{q,me}(\Omega):=
\mathcal E_{q,me,L^\infty_{loc}}(\Omega), \ \ 
\mathcal{N}_{q,me}(\Omega):=
\mathcal N_{q,me,L^\infty_{loc}}(\Omega)\ \ \mbox{ and }\ \ 
$$$$\mathcal{G}_{q}(\Omega):=
\mathcal{E}_{q,me}(\Omega)/\mathcal{N}_{q,me}(\Omega).
$$

We shall call the elements of $\mathcal{E}_{q,me}(\Omega)$ nets of smooth functions with $L^{q}$-moderate growth, while the ones of $\mathcal{N}_{q}(\Omega)$ will be refereed as $L^{q}$-negligible nets. In the same way we define $\mathcal{E}_{q,sm}(\Omega),
\mathcal{N}_{q,sm}(\Omega), \mbox{ and } 
\mathcal{G}_{q,sm}(\Omega)$ as well as the spaces with the continuous representatives with respect to $\varepsilon$ (with subindex $co$). Moreover, we have shown:

\begin{proposition}\label{rel-s-m}
Let $q\in [1,\infty)$. Every $f\in \mathcal G_q(\Omega)$ has a 
representative $(f_\varepsilon)_\varepsilon$ for which the function $(x,\varepsilon)\mapsto f_\varepsilon(x)\in C^\infty(\Omega\times (0,1))$.
\end{proposition}
Thus, we have
\begin{equation}\label{equiv}
\mathcal{G}_{q}(\Omega)=\mathcal{G}_{q,me}(\Omega)=\mathcal{G}_{q,co}(\Omega)=\mathcal{G}_{q,sm}(\Omega), q\in[1,\infty).
\end{equation}
Hence we may always use nets  which are smooth with respect to $\varepsilon$.

%Recall,  $\mathcal{E}_{\infty}(\Omega)=\mathcal{E}_{M,co}(\Omega)$, $\mathcal{N}_{\infty}(\Omega)=\mathcal{N}_{M,co}(\Omega)$ and 
%$\mathcal{G}_{\infty}(\Omega)=\mathcal{G}_{co}(\Omega)$
%are differential algebras that were already considered in the preceding subsection. 
%The reason of this (apparently) exceptional assumption for $q=\infty$ will be clear below. 

When $q=\infty$ , we can define two different spaces associated to the $q=\infty$ (cf. Proposition \ref{spec}):
$$
\mathcal G_{\infty,me}(\Omega)=\mathcal G_{me}(\Omega)  \ \  \mbox{and } \ \  \mathcal {G}_{\infty,co}(\Omega)=\mathcal G_{co}(\Omega)=\mathcal{G}_{\infty,sm}(\Omega)=\mathcal G_{sm}(\Omega)
$$
but Proposition \ref{rel-s-m} does not hold for these two spaces. We shall therefore make a choice for the index $q=\infty$. Our convention is:
$$
\mathcal G_{\infty}(\Omega):=\mathcal G_{co}(\Omega)=\mathcal G_{sm}(\Omega).
$$
Summarizing, Proposition holds for $\mathcal G_{p}(\Omega)$ for every $p\in[1,\infty]$. Furthermore, without lost of generality, we will assume in the sequel that: 

\begin{center}
\emph{all the representatives are continuous with respect to $\varepsilon,$, $q\in[1,\infty].$ }
\end{center}
We also write from now on
$$\mathcal E_q(\Omega):=\mathcal E_{q,co,L^p_{loc}}(\Omega),\;\;
\mathcal N_{q}(\Omega):=\mathcal N_{q,co,L^p_{loc}}(\Omega), \  \  p\in[1,\infty].$$

It is worth mentioning that $\mathcal{D}'(\Omega)$ is embedded into each of the spaces $\mathcal{G}_{q}(\Omega)$ in the same way that it is embedded into the Colombeau algebra of generalized functions. % (\textbf{add more precise comments later}).

Note also that $\mathcal{E}_{q'}(\Omega)\subsetneq \mathcal{E}_{q}(\Omega)$ and $\mathcal{N}_{q'}(\Omega)\subsetneq \mathcal{N}_{q}(\Omega)$ if $q'>q$. %the latter being also true if the subscripts $co$ and $sm$ are considered.
These assertions are shown by Example \ref{exincl1} below. This implies that there exists a canonical linear mapping $\mathcal{G}_{q'}(\Omega)\to\mathcal{G}_{q}(\Omega)$, $q'>q$. As a matter of fact, this mapping is not injective, as the next elementary example shows.

\begin{example}\label{ex1} Consider a net given by $f_{\varepsilon}(x)=n^{-2}e^{n/q}, x\in\mathbb R$ if $\varepsilon\in [n^{-1}-e^{-n},n^{-1}+e^{-n}]$ and $n\geq4$ and $f_{\varepsilon}(x)=0$ otherwise. Then, $(f_{\varepsilon})\in \mathcal{E}_{q}(\mathbb R), q'\leq q$ but $(f_{\varepsilon})\notin \mathcal{E}_{q'}(\mathbb R)$ if $q'>q$. 
%Make obvious changes for continuos and smooth dependence, similar argument for $\mathcal{N}$ and the other assertion of non-injectivity.  
%Details to be added ...
\label{exincl1}
\end{example}

%\subsection{Algebraic properties of the spaces $\mathcal{G}_{q}(\Omega)$} 
%It is well known that $\mathcal{G}_{\infty}=\mathcal{G}_{co}(\Omega)$ is a differential algebra under pointwise multiplication of representatives. 
 The space $\mathcal{G}_{q}(\Omega), 1\leq q<\infty,$ is not an algebra. This follows from Example \ref{ex1} since $(f^2_\varepsilon)_\varepsilon$ does not belong to 
 $\mathcal E_q(\mathbb R).$
Nevertheless, pointwise multiplication on the representative induces a well defined mapping on the corresponding quotients, which operates according to 
$$
(f,g)\in \mathcal{G}_{q'}(\Omega)\times\mathcal{G}_{q}(\Omega)\mapsto f\cdot g\in\mathcal{G}_{r}(\Omega), \ \ \ \frac{1}{q'}+\frac{1}{q}=\frac{1}{r},
$$
as a consequence of H\"{o}lder's inequality. In particular, we obtain the ensuing result.

\begin{theorem} Let $q\in[1,\infty]$. The space $\mathcal{G}_{p}(\Omega)$ is a module over the algebra $\mathcal{G}_{\infty}(\Omega)(=\mathcal{G}_{co}(\Omega))$ under the natural multiplication. Furthermore, it is a differential module, i.e.,
$$
(f\cdot g)^{(\alpha)}=\sum_{\beta\leq \alpha}\binom{\alpha}{\beta} f^{(\beta)}\cdot g^{(\alpha-\beta)}.
$$
\end{theorem}

%\subsection{New classes of generalized constants}
Analogously to the generalized numbers $\tilde{\mathbb{R}}$ and $\tilde{\mathbb{C}}$, one can define the new sets of generalized numbers
$\tilde{\mathbb{R}}_{q}$ and $\tilde{\mathbb{C}}_{q}$, $q\in[1,\infty].$ 
Denote by  $\mathcal E_{0,q}$, resp., $\mathcal N_{0,q}$ spaces of continuous with respect to $\varepsilon$ nets, $(r_\varepsilon)_\varepsilon\in\mathbb C^{(0,1]}$ with the property
$$(\exists a\in\mathbb R)(\int_0^1\varepsilon^{aq}|r_\varepsilon|^q\frac{d\varepsilon}{\varepsilon}<\infty)
$$
resp.,
$$(\forall b<0)(\int_0^1\varepsilon^{bq}|r_\varepsilon|^q\frac{d\varepsilon}{\varepsilon}<\infty).
$$
Then $\tilde{\mathbb{C}}_{q}=\mathcal E_{0,q}/\mathcal N_{0,q}$; $\tilde{\mathbb{R}}_{q}$ is defined with the real nets above.

The sets $\tilde{\mathbb{R}}_{co}=\tilde{\mathbb{R}}_{\infty}$ and $\tilde{\mathbb{C}}_{co}=\tilde{\mathbb{C}}_{\infty}$ are rings, and $\tilde{\mathbb{R}}_{q}$ and $\tilde{\mathbb{C}}_{q}$ become modules over them, respectively.

In particular, $\mathcal{G}_{q}(\Omega)$ becomes a module over the ring of generalized constants $\tilde{\mathbb{C}}_{co}$.

\section{Besov type spaces of generalized functions \cite{psvv}}\label{bes}

In the sequel we will consider representatives of generalized functions consisting of continuous functions with respect to $\varepsilon\in(0,1],$ as discussed in the previous section.

Let $p\in[1,\infty], q\in[1,\infty)$. We consider $(f_\varepsilon)_\varepsilon\in{\mathcal E}_{q}(\Omega)
$ 
such that for given $k\in \mathbb {N}$ and $s\in\mathbb{R}$ there holds 
\begin{equation}
\label{eqnetg1}
(\forall \omega\subset \subset \Omega)(\int_0^1\varepsilon^{sq}||f_\varepsilon||_{W^{k,p}(\omega)}^qd\varepsilon/\varepsilon<\infty).
\end{equation}

We say that a net $(f_{\varepsilon})_{\varepsilon}\in \mathcal{E}_{q}(\Omega)$ belongs to $\mathcal{E}^{k,-s}_{q,L^{p}_{loc}}(\Omega)$ if (\ref{eqnetg1}) holds. Special attention will be devoted in Section \ref{zcass} to Zygmund type spaces in the case $q=\infty$.

 Furthermore,  for $k=\infty$, we put  $$\mathcal{E}^{\infty,-s}_{q,L^{p}_{loc}}(\Omega)=\bigcap_{k\in\mathbb{N}}
\mathcal{E}^{k,-s}_{q,L^{p}_{loc}}(\Omega).$$ %When $p=\infty$, $q<\infty$, we set $\mathcal{E}^{k,-s}_{q,loc}(\Omega)=\mathcal{E}^{k,-s}_{q,L^{\infty}_{loc}}(\Omega)$.
The spaces of the following definition will be vital in our study of Besov type regularity.

%\marginpar{Part of the definition was moved to the text, it was extremely heavy to read.}  
\begin{definition} \label{def1} Let $s\in\mathbb{R}$, $k\in\mathbb{N}_{0}\cup\left\{\infty\right\}$, $q\in [1,\infty]$, and $p\in [1,\infty]$. Then $\mathcal{G}^{k,-s}_{q,L^{p}_{loc}}(\Omega)$ is  the quotient space 
\begin{equation}\label{defeq1}
\mathcal{G}^{k,-s}_{q,L^{p}_{loc}}(\Omega)=
 \mathcal{E}^{k,-s}_{q,L^{p}_{loc}}(\Omega)/\mathcal{N}_{q}(\Omega).
\end{equation}
\end{definition}

We have  ${\mathcal{G}}^{k,-s}_{q,L^{p}_{loc}}(\Omega)\subset \mathcal{G}_{q}(\Omega)$ for any $p\in[1,\infty]$. Note that the definition does not depend on the representatives.

We list some properties of these vector spaces of generalized functions in the next proposition. %Their proofs follow immediately from Definition \ref{def1} and our previous results.

\begin{proposition}
\label{propc1} Let $s\in\mathbb{R}$, $k\in\mathbb{N}_{0}\cup \left\{\infty\right\}$, $q\in[1,\infty]$ and $p\in[1,\infty]$. 
\begin{itemize}
\item [(i)] $
{\mathcal{G}}_{q,L^p_{loc}}^{k,-s}(\Omega)\subseteq
{\mathcal{G}}_{q,L^p_{loc}}^{k_1,-s_1}(\Omega)$ if and only if $k\geq k_1$ and $s\leq
s_1$.
\item [(ii)] Let $P(D)$ be a
differential operator of order $m\leq k$ with constant coefficients. Then
$P(D):{\mathcal{G}}_{q,L^{p}_{loc}}^{k,-s}(\Omega)\rightarrow
{\mathcal{G}}_{q,L^{p}_{loc}}^{k-m,-s}(\Omega).$
\item[(iii)] Let $ \infty>r>q,  \rho<s.$ Then $\mathcal G^{m,-s}_{q,L^p_{loc}}(\Omega)\subset \mathcal G^{m,-\rho}_{r,L^p_{loc}}(\Omega).$ 
\item[(iv)] $\mathcal{G}^{\infty,-s}_{q,L^{p}_{loc}}(\Omega)=\mathcal{G}^{\infty,-s}_{q,L^{\infty}_{loc}}(\Omega),
$
\end{itemize}
\end{proposition}

\section{Characterization of  Besov regularity of distributions \cite{psvv}}
\label{cass}

We defined in \cite{p-r-v} non-degenerate wavelets and generalized Littlewood-Paley (LP) pairs of order $\alpha\in\mathbb R.$ Here we will simplify
the exposition considering  special (LP) pairs defined as follows. A (LP) pair is $(\phi_1,\psi_1),$ where
$\hat{\phi}_1, \hat{\psi}_1\in\mathcal{D}(\mathbb R^d),$ 
$ \hat \phi_1\equiv 1 $ in a ball $B(0,r)$, 
$ \mbox{ supp }\hat \psi_1 \subset B(0,r_1)\setminus B(0,r_2), r_1>r>r_2, $ 
$\hat \psi_1\equiv 1$ in a neighborhood of $S(0,r)$
($B(0,\rho)$  denotes an open ball whose boundary is the sphere $S(0,\rho)$).

Clearly this definition is satisfied by the special (LP) pair $(\varphi,\psi)$ from Section \ref{gfc}.
In the sequel, we will assume that the mollifier $\phi$ is chosen so that 
\begin{equation}\label{lp}
(\varphi*\phi,\psi*\phi)
\end{equation}
makes an (LP) pair.

\begin{proposition}\label{est}

Let $(\phi_1,\psi_1)$ be a Littlewood-Paley pair as above.
%$\hat{\phi}_1\in C^\infty_0(\mathbb R^n)$ and  $\hat{\psi}_1\in C^\infty_0(\mathbf R^n)\setminus\{0\}$ be a compatible pair. 
Define  $||u||^s_{(1),q,p}$ with this pair,  see (\ref{besov}).
Then, for $k=|\alpha|,$
$$||\psi_\varepsilon*u^{(\alpha)}||_{L^p(\mathbb R^d)}
\leq C||(\psi_1)_\varepsilon *u^{(\alpha)}||_{L^p(\mathbb R^d)}.$$
In particular the norms $||u||_{q,p}^{s}$ and 
 $||u||_{(1),q,p}^{s}$ are equivalent.
\end{proposition}

The space $B^{s}_{q,L^{p}}(\mathbb R^d)\cap\mathcal{E}'(\Omega)$ is 
naturally embedded into $\mathcal{G}_{q}(\Omega)$, through convolution 
with a mollifier: $T\mapsto T_\varepsilon=T\ast \phi_{\varepsilon}|_{\Omega}, {\varepsilon}\in(0,1)$. 

\begin{proposition} \label{converse}
Let $f\in\mathcal E'(\mathbb R^d)$ 
such that  $$\iota(f)\in\mathcal \iota(\mathcal D'(\mathbb R^d))\cap \mathcal G^{0,s}_{q,L^p_{loc}}(\mathbb R^d).
$$
Then $f\in B^{s}_{q,L^p}$.  
\end{proposition}

\begin{theorem} Let $s>0$ and $k\in\mathbb{N}$. Then,
$$\mathcal{G}^{k,-s}_{q,L^{p}_{loc}}(\mathbb R^d)\cap \iota(\mathcal{E}'(\mathbb R^d))\supset\iota (B^{k-s_0}_{q,L^{p}})$$
for any $s_0<s.$
\end{theorem}

Recall
% (\cite{col1}, \cite{ober001})
 that a net $(f_{\varepsilon})_{\varepsilon}\in\mathcal{E}^{(0,1)}(\Omega)$, or the generalized function
  $f=[(f_{\varepsilon})_{\varepsilon}]$, is strongly  associated  to  $T\in\mathcal D'(\Omega)$ if there exists  $b>0$ such that
   %$\lim_{\varepsilon\to0}f_{\varepsilon}=T$ in the weak topology of $\mathcal{D}'(\mathbb R^d)$, that is,
\begin{equation}
\label{eqnet1}
(\forall \rho \in \mathcal{D}(\Omega))(\langle T-f_\varepsilon, \rho
\rangle=o(\varepsilon^b),\ \varepsilon\to 0).
\end{equation}
%\marginpar{I corrected and made more general definitions and statements, many improvements are %still needed. Maybe you %should not present the results of this subsection yet in the conference}

With $o(1)$ instead $o(\varepsilon^b)$ in (\ref{eqnet1}), one has the notion of weak association.
We introduce a new  concept of  association.
\begin{definition}\label{qaa}
Let $T\in\mathcal D'(\mathbb R^d)$, $(f_\varepsilon)_\varepsilon\in \mathcal E_{q}(\mathbb R^d)$. We say that $(f_\varepsilon)_\varepsilon$ is strongly $ q-$associated to $T$ if there exists $b>0$ such that
\begin{equation}\label{ass}
(\forall \rho \in \mathcal{D}(\mathbb R^d))(\int_0^1\varepsilon^{-bq}|\langle T-f_\varepsilon, \rho
\rangle|^qd\varepsilon/\varepsilon<\infty).
\end{equation}
%for every $\omega\subset\subset \Omega$
\end{definition}

Clearly the strong association implies the $q$-association and the converse does not hold. 
Moreover the weak association and the $q$ associations are not comparable.

\begin{theorem}
\label{not1}
Let $T\in\mathcal E'(\mathbb R^d)$ and $[(f_\varepsilon)_\varepsilon]\in \mathcal G^{k,s}_{q,\infty}(\mathbb R^d)$
 for some $k\in\mathbb N$ and every  $s>0.$
Assume that $T $ and $(f_\varepsilon)_\varepsilon$ are strongly $q$-associated, $q\geq 1$.
Then $\iota(T)\in\mathcal{G}_{q,L^{\infty}_{loc}}^{k,s}$ for  every  $s>0$ . In particular, $T\in B^{k+s}_{q,\infty}$ for every $s>0.$
\end{theorem}

\section{Zygmund  regularity through association \cite{psvmos}}
\label{zcass}
In this section we present the local regularity of distributions in connection with the Zygmund type classes $\mathcal{G}^{k,-s}_{\infty,L^{\infty}_{loc}}(\Omega)$  (see (7.1) with $q=\infty$) . 

The next  theorem provides a precise  characterization of those distributions that 
belong to $\mathcal{G}^{k,-s}_{\infty,L^\infty_{loc}}(\Omega)$,
 denoted in the sequel as $\mathcal{G}^{k,-s}(\Omega)$, they turn out to be elements of a Zygmund space. We only consider the case $s>0$, since for $s\leq 0,$ one has $\mathcal{G}^{k,-s}(\Omega)\cap\iota(\mathcal{D}'(\Omega))=\left\{0\right\}$.

\begin{theorem}
\label{1mtheorem} Let $s>0$. We have $\mathcal{G}^{k,-s}(\Omega)\cap\iota(\mathcal{D}'(\Omega)) =\iota(C_{\ast,loc}^{k-s}(\Omega))$. 
\end{theorem}

 This implies that for $r\in\mathbb{R}$ and any non-negative integer $k>r$,
 $$\iota(C_{\ast,loc}^{r}(\Omega))=\mathcal{G}^{k,r-k}(\Omega)\cap\iota(\mathcal{D}'(\Omega)).$$

Consequently,  we immediately have $\iota(\mathcal{D}'(\Omega))\cap \mathcal{G}^{\infty}(\Omega)=\iota(C^{\infty}(\Omega))$.

We return to the strong association but with  the more general rate of approximation in (\ref{eqnet1}).
 Let $R:(0,1]\to \mathbb{R}_{+}$ be a positive function such that $R(\varepsilon)=o(1),$ $\varepsilon\to0.$ We write $
T-f_{\varepsilon}=O(R(\varepsilon)) \ \mbox{  in } \mathcal{D}'(\Omega)
$
if 
\begin{equation*}
(\forall \rho \in \mathcal{D}(\Omega))(\langle T-f_\varepsilon, \rho
\rangle=O(R(\varepsilon)),\ \varepsilon<1).
\end{equation*}

We now present our results concerning the regularity analysis through association. 

\begin{theorem} \label{regas}
 Let $T\in\mathcal D'(\Omega)$ and let $f=[(f_\varepsilon)_\varepsilon]\in\mathcal{G}({\Omega})$ be associated to it. Assume that $f\in\mathcal{G}^{\infty}(\Omega)$. If 
$(f_{\varepsilon})_{\varepsilon}$ approximates $T$ with convergence rate:
\begin{equation}
\label{req3}
(\exists b>0)(T-f_{\varepsilon}=O(\varepsilon^{b}) \mbox{ in }\mathcal{D}'(\Omega)).
\end{equation} 
Then $T\in C^\infty(\Omega).$
\end{theorem}

 \begin{theorem}\label{mtheorem}
 Let $T\in\mathcal{D}'(\Omega)$ and let $f=[(f_\varepsilon)_\varepsilon]\in\mathcal{G}(\Omega)$ be a net of smooth functions associated to it. Furthermore, let $k\in\mathbb{N}$. Assume that either of following pair of conditions hold:
\begin{enumerate}
\item[\textnormal{(i)}] $f\in \mathcal{G}^{k,-a}(\Omega)$, $\forall a>0$, namely,
\begin{equation}
\label{req2}
(\forall a>0)(\forall \omega\subset \subset \Omega)(\forall \alpha \in {\mathbb
N}^d,|\alpha|\leq k)(\sup_{x\in\omega}|f^{(\alpha)}_\varepsilon(x)|=O(\varepsilon^{-a})),
\end{equation}
and the convergence rate of  $(f_{\varepsilon})_{\varepsilon}$ to $T$ is as in (\ref{req3}).
 
\item [\textnormal{(ii)}] $f\in \mathcal{G}^{k,-s}(\Omega)$ for some $s>0$, and there is a rapidly decreasing function $R:(0,1]\to\mathbb{R}_{+}$, i.e., $(\forall a>0)(\lim_{\varepsilon\to0}\varepsilon^{-a}R(\varepsilon)=0)$, such that
\begin{equation}
\label{req4}
T-f_{\varepsilon}=O(R(\varepsilon))\  \mbox{ in }\mathcal{D}'(\Omega).
\end{equation}
\end{enumerate}
Then, $ T\in C_{*,\:loc}^{k-\eta}(\Omega)$ for every $\eta>0$. 

\end{theorem}

\end{document}